\documentclass[%
reprint,
showpacs,preprintnumbers,
 amsmath,amssymb,
 aps,
pra,
]{revtex4-1}

\usepackage[usenames]{color}
\usepackage{fancyhdr}
\usepackage[normalem]{ulem}

\usepackage{graphicx}
\usepackage{dcolumn}
\usepackage{bm}

\begin{document}


\title{Detecting the temporal structure of intermittent phase locking}


\author{Sungwoo Ahn}
\author{Choongseok Park}
\affiliation{
Department of Mathematical Sciences and Center for Mathematical Biosciences, Indiana University Purdue University Indianapolis, IN 46032, USA \\
}

\author{Leonid L. Rubchinsky}
\email{leo@math.iupui.edu}
\affiliation{
Department of Mathematical Sciences and Center for Mathematical Biosciences, Indiana University Purdue University Indianapolis, IN 46032, USA \\
Stark Neurosciences Research Institute, Indiana University School of Medicine, Indianapolis, IN 46032, USA \\
 }


\begin{abstract}

This study explores a method to characterize temporal structure of intermittent phase locking in oscillatory systems. When an oscillatory system is in a weakly synchronized regime away from a synchronization threshold, it spends most of the time in parts of its phase space away from synchronization state. Therefore characteristics of dynamics near this state (such as its stability properties/Lyapunov exponents or distributions of the durations of synchronized episodes) do not describe system's dynamics for most of the time. We consider an approach to characterize the system dynamics in this case, by exploring the relationship between the phases on each cycle of oscillations. If some overall level of phase locking is present, one can quantify when and for how long phase locking is lost, and how the system returns back to the phase-locked state. We consider several examples to illustrate this approach: coupled skewed tent maps, which stability can be evaluated analytically, coupled R\"{o}ssler and Lorenz oscillators, undergoing through different intermittencies on the way to phase synchronization, and a more complex example of coupled neurons.  We show that the obtained measures can describe the differences in the dynamics and temporal structure of synchronization/desynchronization events for the systems with similar overall level of phase locking and similar stability of synchronized state.
\end{abstract}


\maketitle


\section{\label{intro}Introduction}

Synchronization of oscillations is a widespread phenomenon in
physical sciences and beyond~\cite{arkady_pikovsky, boccaletti2002},
including physiology~\cite{glass2001} and
neuroscience~\cite{engel2001, buzsaki2004}, and has being studied
using approaches and methods of physics and nonlinear dynamics even
in that applied context~\cite{nowotny2008, rabinovich2006}. In
particular, phase synchronization, dynamical phenomena where the
phases of oscillations approach each other with time while both
phases may remain chaotic and their amplitudes remain
uncorrelated, appears to be quite common in nature and is
well-studied~\cite{arkady_pikovsky}. It is known that for some
parameter values (usually, moderate coupling strength), phase
synchronization can occur in an intermittent fashion~\cite{lee_kwak_lim1, ring, rim_kim_kang}. It will present itself as an intermittent phase locking.
The phase difference of two
oscillators is close to constant during finite time
intervals. These intervals of temporal synchronization epochs are
random and interspersed by desynchronization events which are
characterized by varying phase differences. Different types of
intermittent behaviors may take place in different
systems~\cite{rim_kim_kang} and transition through intermittency may
depend on the geometric structure of chaotic
attractors~\cite{zhao2005}. Three types of intermittencies have been
observed near the onset of phase synchronization: the type-I
intermittency~\cite{lee_kwak_lim1}, the eyelet
intermittency~\cite{lee_kwak_lim1, uni_eyelet}, and the ring
intermittency~\cite{ring}.

Different approaches to synchronization analysis exist. One may
approach the problem from the time-series analysis perspective,
characterizing how much time-series are correlated in some
particular sense. This view essentially tells how close we are to
the synchronized state. Another approach is to study the stability
of the synchronized state, via, for example, Lyapunov exponents, so
that small positive Lyapunov exponents will describe the weak
instability and, potentially, intermittency of synchronization.
Lyapunov exponents, distribution of durations of synchronized
episodes etc. are ways to explore the synchronized
state/synchronization manifold. Certain universality of behavior of
a weekly unstable synchronization is expected and is captured by
these approaches. Studies of Lyapunov exponents spectra ~\cite{zero_lyapunov} may go further and explore the synchronous/desynchronous states dynamics in the vicinity of the saddle-node bifurcation point.
However, these approaches become less powerful when the system is far away from bifurcation point and do not describe the
desynchronized episodes of intermittent or otherwise variable and weak synchronization in a wide range of parameter space. The desynchronized dynamics (``reinjection mechanisms" in intermittency terminology) may be not universal and
depend on the peculiarities of the system under study.

If the system is in imperfectly synchronized state, but is still close to synchronization threshold (the case, which is probably expected in many physics or engineering applications), it spends most of the time in an almost synchronized state (near synchronization manifold) and the focus on this state is natural. However, the synchrony can be weak. In particular, in living systems the synchronization may be very imperfect and weak, and too strong synchrony (high degree of temporal coordination) may be responsible for pathological states, such as schizophrenia or Parkinson's disease~\cite{schnitzler2005, uhlhaas2006}. In terms of the phase space, the system spends relatively small fraction of time near synchronization manifold in the phase space. Thus even if Lyapunov exponents or distributions of the durations of synchronized episodes may be informative about this part of dynamics, they do not tell much about dynamics overall, because it is mostly desynchronized. This calls for the study of the highly-variable synchrony and for the study of the structure of the phase space away from the relatively strongly unstable synchronization manifold. One may not expect much of universality here, but the characterization of this dynamics should be possible.

One approach is to characterize the degree of synchronization for a series of sliding time-windows  and obtain statistical estimates of significance (as in, for example,~\cite{hurtado_rubchinsky_2004}). This may reveal important details of synchronized dynamics, not accessible otherwise (as in~\cite{hurtado2004} for the case of tremor oscillations). However, synchronization is not instantaneous phenomena, so this approach is not designed to detect changes in the dynamics on the very short time-scales. Yet, changes of this kind (very intermittent, variable synchronization) have been observed experimentally and have been conjectured to have functional significance~\cite{park_rubchinsky1, park_rubchinsky2}.

Another example, where the fine temporal details of synchrony may be relevant may come from ecological dynamics. Synchronization in ecological dynamics is hardly perfect. It was conjectured that prolong zero-lag synchronization elevates risks of extinction, so that there is an interest in these cases (reviewed in, e.g. ~\cite{ecolo-synch}). One may want to distinguish between the moderate level of synchrony with few prolong synchronized episodes and many short synchronized episodes.

Here we expand on the approach of ~\cite{park_rubchinsky1} and study a method to analyze intermittent and weak phase locking from the nonlinear dynamics perspective. This approach is based on the constructions of the first-return maps for the difference of phases, partition of the resulting phase space into synchronized and non-synchronized regions and the analysis of transitions between them. We show that for the case when some overall level of phase locking is present, one can characterize how this phase locking (properly defined) is gained, maintained or lost for each cycle of oscillations and what it means in terms of the organization for the phase space of coupled oscillators system. We do so with examples drawn from several typical model systems, ranging from simple coupled tent maps to more realistic conductance-based models of coupled cells. We also show how one can characterize the differences in the dynamics of coupled systems, which undergo through the same kind of intermittency and posses the same stability properties of synchronization manifold.

The structure of this paper is the following. Section~\ref{methods_sync} presents the set up of the method. Section~\ref{coupled_maps_example} considers the method as applied to the dynamics of simple coupled maps. Section~\ref{intermittency} considers different types of intermittencies for phase synchronization. Section~\ref{neuronal_model} illustrates the ideas of the paper for a more complex dynamical system. Finally, in Section~\ref{discussion} we discuss the main ideas behind the method, summarize its properties, and discuss its applicability and possible further developments.

\section{\label{methods_sync} Methods to study the temporal variability of the synchronization}

Suppose that two oscillatory signals allow for a reasonably good reconstruction of their phases $\phi_{1}(t)$ and $\phi_{2} (t)$. In a discrete-time version (and thus applicable to experimentally obtained data) one can consider a standard index to characterize the strength of the phase locking between these two signals:
\begin{equation}\label{sync_gamma} \gamma = || \frac{1}{N} \sum_{j=1}^N e^{i \theta (t_j)} ||^2,\end{equation}
where $\theta (t_j) = \phi_{1}(t_j) - \phi_{2} (t_j)$ is the phase difference. This synchronization index varies from 0 to 1, that is from complete lack of phase locking to perfect phase locking (and does not tell anything about amplitudes of the signals). One may also compute Lyapunov exponents, which characterize stability of the synchronization manifold. In general, negative values of all transversal Lyapunov exponents are not necessarily sufficient to guarantee perfect synchronization~\cite{arkady_pikovsky}, nevertheless their values describe the degree of stability/instability of the synchronized state and the strength of synchronization.

As we explained in the Introduction, these measures are not designed to describe the fine temporal structure of the dynamics, which motivates the following approach (its major steps were sketched by us in~\cite{park_rubchinsky1} and its explanation, validation and properties are presented in subsequent sections of this paper). First, one needs to confirm that there exist some degree of phase locking between the phases of two variables $x$ and $y$. This can be detected, for example, by computing the Eq.~(\ref{sync_gamma}) and comparing the resulting value of $\gamma$ against some kind of significance value, computed analytically or via surrogates~\cite{hurtado_rubchinsky_2004}.
It should be emphasized here that our method allows for an analysis of the temporal development of phase difference if some level of synchrony (some preferred phase-locking angle) is present. It makes no sense if there were no
synchrony between the signals at all.

Then set up a check point for the phase of $y$: whenever the phase
of $y$ crosses this check point from below, record the phase value
of $x$, generating a set of consecutive phase values $\{\phi_{x,
i}\},$ $i=1, \ldots, N$, where $N$ is the number of such level
crossings in an episode under analysis. Then plot $\phi_{x, i+1}$
{\it{vs.}} $\phi_{x, i}$ for $i=1, \ldots, N-1$. The properties of
these plots will characterize the dynamics of the synchronization. A
tendency for predominantly synchronous dynamics will appear as a
cluster of points, with the center at the diagonal $\phi_{x, i+1} =
\phi_{x, i}.$

\begin{figure}[htp]
   \centering
   \includegraphics[width=3.4in,height=2.7in]{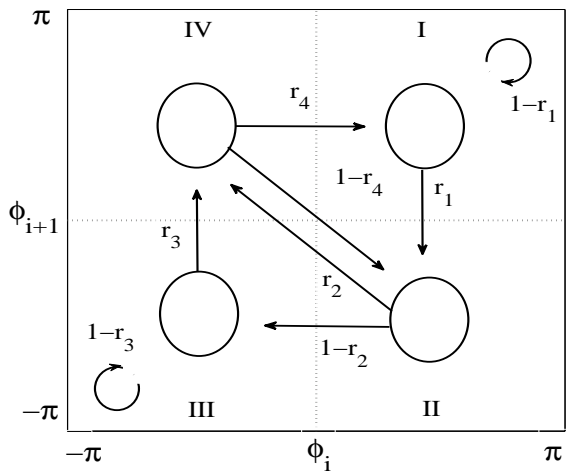}
  \caption{\label{phase_diagram} Diagram of the $(\phi_{x, i}, \phi_{x, i+1})$ first-return map. The arrows indicate all possible transitions from one region to another and the expressions next to the arrows indicate the rates for these transitions. }
\end{figure}

After determining the center of the cluster (that is, determining the preferred phase difference between two signals, at which they have the tendency to be locked) all values of the
phases are shifted to a position with the center of the first region
(quadrant) (see Fig.~\ref{phase_diagram}). This step is not necessary and is done here for the uniformity of the analysis.
Then we consider how the system
leaves the synchronization cluster and its vicinity and how it
returns back to the synchronization by quantifying transitions
between different regions of the $(\phi_{x, i}, \phi_{x, i+1})$
space.

\begin{figure}[htp]
   \centering
   \includegraphics[width=3.4in,height=3.4in]{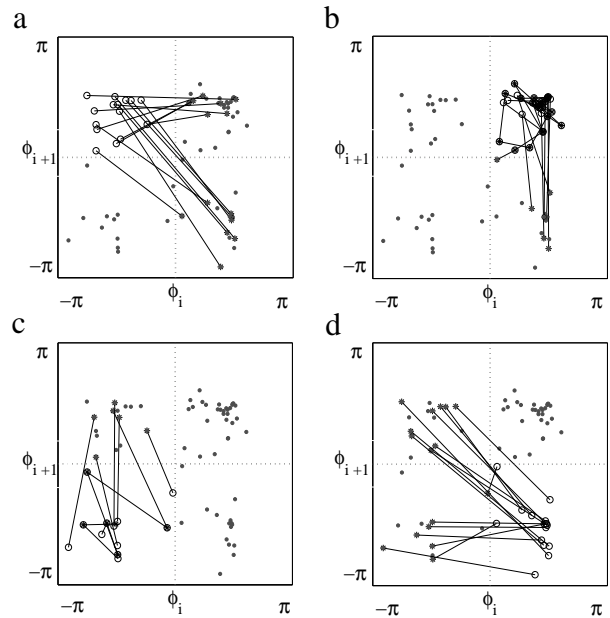}
  \caption{\label{return_map_example} An example of the first-return map for two coupled Lorenz oscillators. All four first-return plots have the same data points (gray circles), but each subplot (a-d) presents the evolution of points from one region (I, II, III, and IV, respectively). If a point evolves from one region to another region, then we represent it as $\circ-\ast$. If a point evolves within the same region, then we represent it as $\circ-\circ$. Thus each plot shows the transitions from a corresponding part of the phase space. We can compute the transition rates $r_{1, 2, 3, 4}$ from the panels (b), (d), (c), and (a), respectively. }
\end{figure}

This phase space is then partitioned into four regions as illustrated in Fig.~\ref{phase_diagram} (let us remind, that the space of $(\phi_{x, i}, \phi_{x, i+1})$ is in fact torus). While the system is in the first region, it is considered to be in a synchronized state. Dynamics outside of the first region will be called a desynchronized state. Thus the synchronized state here is the one, where the deviation from the preferred phase angle is less than $\pi/2$. This is a compromise value: not too small to allow for some moderate fluctuations of the phases; not too large to allow for the phases to be sufficiently related and to be involved in some function. It also allows for symmetric partitioning of the phase space and definition of a few easy-to-compute transition rates (see below).
Some situations may require a different definition of the synchronized state and the rest of the method will need to be modified. Nevertheless, even though this partition is quite coarse, it allows for an inspection of
the fine temporal details, as we will show below. Four resulting regions in the phase space are numbered in a clockwise manner (Fig.~\ref{phase_diagram}), since this is the primary direction of the dynamics. An example of the first-return map for two coupled Lorenz oscillators is presented in Fig.~\ref{return_map_example} (detailed study of this system with the method described here is in the Section~\ref{intermittency}, see Fig.~\ref{Lorenz_Type_I_elelet_best1}). 

Transition rates $r_i$ ($i =1, 2,
3, 4$) for the transitions between four regions of the map  are defined as the number
of points in a region, from which the system leaves the region,
divided by the total number of points in that region. For example,
$r_1$ is the ratio of the number of trajectories escaping the first
region for the second region to the number of all points in the first
region. While the time-averaged measures of synchrony characterize whether the
synchronization is strong or weak overall, utilization of these
rates lets us explore the dynamics of the synchronization in time.

Further exploring the properties of the dynamics in the space $(\phi_{x, i}, \phi_{x, i+1})$, one may compute the relative frequencies of desynchronization events of different durations. In the considered first-return map approach, the duration of a desynchronization event is the number of steps that system spends away from the first region minus one (because, the point on the map has two coordinates, one of which is a ``future" phase). The shortest duration of the desynchronization event corresponds to the shortest path 2-4-1 (note that the desynchronization will always start at the second region by the virtue of our description of the dynamics). This will correspond to the desynchronization length of one cycle (in other words, in two cycles the phases are back in a locked state). If we want to consider the chance of the duration of two cycles, we should consider the path 2-3-4-1. Longer desynchronization events will have many different paths corresponding to them. We will systematically compute the relative frequencies of the desynchronization events of different durations in the Section~\ref{intermittency} for different systems.

The transition rates $r_{2, 3, 4}$ are related to the durations of desynchronizations. Higher values of $r_{2, 3, 4}$ promote short desynchronization episodes, while their low values promote long desynchronization events (note that for a well-synchronized dynamics, there may be very few points in the third quadrant, and numerical estimation of $r_3$ may be unreliable).
If transitions are independent then the relative frequencies of the durations of desynchronization events can be estimated by multiplying the corresponding rates of the transitions for different paths. For example, to estimate
the probability of a desynchronization event of the shortest
duration we should consider the shortest path 2-4-1 and the
corresponding probability is $r_2 r_4.$ This will correspond to the
desynchronization length of one cycle. Longer desynchronizations will require summation of the products of rates corresponding to longer paths of equal length. In general, this independence (or near independence) may or may not be in place, but we consider several examples, where the difference is not large and the rates $r_{2, 3, 4}$ really describe the desynchronization durations.

The rate $r_1$ is related to the averaged duration of the laminar
phase and thus characterizes the property of synchronized state. The
average duration of the synchronized phase ~$<\ell>$ (for an
intermittency near onset of phase synchronization, this quantity
usually scales in a universal manner for specific intermittency
types) is inversely proportional to $r_1$:
\begin{equation} <\ell>  \sim 1/r_1 \end{equation}
and $<\ell> = 1/r_1$, if $<\ell>$ is measured in the number of
iterations of the map $\phi_{x, i}$, which is essentially the number
of cycles of oscillations. The closer $r_1$ is to 1, the more
frequently the synchronized dynamics is interrupted.

Thus the rates
permit evaluation of whether the weakness of synchrony is due to a
few relatively long desynchronization events or to large number of
short desynchronization events.

\section{\label{coupled_maps_example} Coupled maps example}

In this section we consider behavior of rates for different parameter values (and thus for different transversal Lyapunov exponents of synchronized state and the synchronization index $\gamma$) in a simple system. We consider linearly coupled skew tent maps. This relatively simple map allows for exact computation of Lyapunov exponents. The properties of the phase space depend on the parameters in a simple way. Thus we are able to see the relationship between system's parameters, dynamics, Lyapunov exponents and the rates $r_i$. This system is not necessarily an appropriate choice for the study of phase synchronization as a dynamical phenomenon, when phases are correlated and amplitude are not~\cite{phase-discrete}. However we are not as much concerned with the phase synchronization per se in this sense, as we are concerned with a simple treatable example of synchronous dynamics here. Therefore, coupled skew tent maps prove a possibility to illustrate our ideas in a simple system, where analytical treatment is partially possible.

Consider a skew tent map
\begin{eqnarray}\label{tent_map_example}
f(a, x)=
\begin{cases}
\frac{x}{a} & \text{ if } 0 \leq x \leq a, \\
\frac{1-x}{1-a} & \text{ if } a < x \leq 1,
\end{cases}
\end{eqnarray}
where $0< a < 1.$ We will use two such maps with
variables $x$ and $y$. Let us consider the following coupling:
\begin{eqnarray}\label{basic_coupled_system} x(t+1) &=& (1-\varepsilon) f(a, x(t)) + \varepsilon f(a, y(t)), \\
y(t+1) &=& \varepsilon f(a, x(t)) + (1-\varepsilon) f(a, y(t)) \nonumber , \end{eqnarray} where $\varepsilon$ is the coupling strength.

To characterize synchrony between $x$ and $y$, we consider a new variable
\begin{equation}\label{theta_diff} \theta(t) =y(t) - x(t),\end{equation}
which will be an analog of the phase difference described above. In these coupled maps, synchronous state $x=y$ is stable for $\varepsilon > \varepsilon_c$ where the critical value $\varepsilon_c < 1/2.$
Note that $\theta$ in Eq.~(\ref{theta_diff}) varies between $-1$ and $1$ while the phase variable $\phi_i$ in Fig.~\ref{phase_diagram} varies between $-\pi$ and $\pi.$ In this section, we rescale the phase diagram accordingly to define the transition rates $r_i.$

Two Lyapunov exponents can be computed analytically (see e.g. \cite{arkady_pikovsky}):
\begin{eqnarray}\label{two_lyap} \lambda (a) &=& -a \ln a - (1-a)
\ln (1-a),
\\ \nonumber
        \lambda_\perp (a, \varepsilon)  &=& -a \ln a - (1-a) \ln (1-a) +
        \ln|1-2\varepsilon|. \nonumber \end{eqnarray}
The first Lyapunov exponent $\lambda (a)$ corresponds to the motion governed by the
one-dimensional dynamical system $x(t+1) = f(a, x(t))$ whose
subspace is $\{(x, y) ~|~ x = y \}.$ The second, transversal Lyapunov exponent
$\lambda_\perp (a, \varepsilon)$ characterizes the dynamics transverse
to that invariant subspace.  Note that $\lambda(a)$ is
always positive for $0 < a <1$ and has a maximum $\ln 2$ at $a=1/2.$
Since  $\lambda_\perp(a, \varepsilon) = \lambda (a) +
\ln|1-2\varepsilon|,$ the transversal Lyapunov exponent $\lambda_\perp(a, \varepsilon)$
can change from positive to negative depending on $a$ and
$\varepsilon.$ Both Lyapunov exponents are symmetric around~$a=1/2.$ That is, $\lambda(a)
= \lambda(1-a)$ and $\lambda_\perp(a, \varepsilon) = \lambda_\perp(1-a,
\varepsilon)$. The synchronization stability threshold $\lambda_\perp(a,
\varepsilon_c)=0$ is reached at $\varepsilon_c = \frac{1}{2} - \frac{1}{2} a^a
(1-a)^{(1-a)}.$

It is interesting to note that the coupled maps (Eqs.~(\ref{tent_map_example}) and~(\ref{basic_coupled_system})) have exact
solutions when the system is at a bifurcation point and is not
robust. For example, for with $a=1/2$ and $\varepsilon=1/4$ (and
thus $\lambda = \ln 2$ and $\lambda_\perp = 0$) initial conditions
$(x,y)=(0.46, 0.5)$ lead to a period-20 trajectory and initial
conditions $(x, y)= (0.044, 0.052)$ lead to a period-100 trajectory.
These trajectories are not hard to compute, and then one can formally compute
$r_i$ ($r_1 = 2/5,$ $r_2=3/5,$ $r_3=2/5,$ and $r_4=2/5$ for the
former trajectory and $r_1 = 13/25,$ $r_2 = 12/25,$ $r_3 = 13/25,$ and
$r_4 = 13/25$ for the latter). However, the system does not exhibit
synchronous behavior, there is no single preferred value for the
 difference of variables $\theta$, and computation of rates is
not appropriate.

\begin{figure}[htp]
    \centering
   \includegraphics[width=3.4in,height=3in]{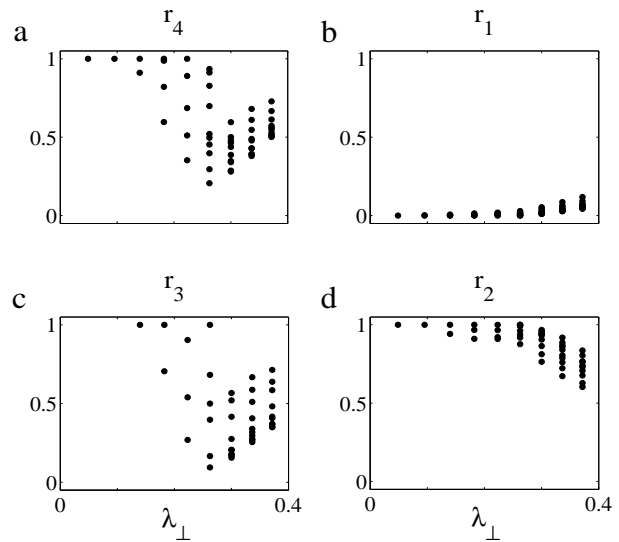}
  \caption{\label{tent_map_rates_trans_lyap1}
We varied parameters $a$ and $\varepsilon$
  to compute the transversal Lyapunov exponent $\lambda_\perp$ and rates $r_i$. The rate $r_1$ ($r_2$)
  and $\lambda_\perp$ have  a positive (negative) relation. Rates $r_3$ and $r_4$ present more spread and more complex relationship with $\lambda_\perp$.}
\end{figure}

We numerically iterate the coupled maps
(Eqs.~(\ref{tent_map_example}) and~(\ref{basic_coupled_system})) and
compute the rates $r_i$ for different $a$ and $\varepsilon.$ We vary $a$
from $0.125$ to $0.35$ with step size $0.025$ and $\varepsilon = (1-
k a^a (1-a)^{1-a})/2$ where $k$ varies from $1.05$ to $1.45$ with
step size $0.05$ (with this set-up $\lambda_\perp = \ln k$).
Fig.~\ref{tent_map_rates_trans_lyap1} shows the relationship between
$\lambda_\perp$ and the rates~$r_i$. Each dot represents a different
pair of $(a, \varepsilon).$ Since $\lambda_\perp$ characterizes the
evolution of the trajectory transverse to the invariant subspace
$\{(x, y) ~|~ x = y\},$ the smaller $\lambda_\perp$ corresponds to
less unstable manifold and more synchronized dynamics, which results
in the smaller~$r_1.$ As $\lambda_\perp$ decreases, the points $(x,
y)$ coalesce to the diagonal line $x=y$ in the $(x, y)$ plane. This
shortens the duration of desynchronization events and increases
the reinjection probability to the synchronous state. Thus, $r_1$
tends to decrease while $r_{2, 4}$ tend to increase as
$\lambda_\perp$ decreases. However, one can clearly see that the same value
of $\lambda_\perp$ may correspond to drastically different rates and
thus to different timing of synchronized/desynchronized dynamics.

If there is no point in the given region, then the transition rate
at that region is undefined. As $\lambda_\perp$ decreases, the
system spends less and less time in desynchronous regions. This
makes the number of points in the third region extremely small (or
sometimes zero). Thus numerically computed values of $r_3$ are less
reliable in this case (or are undefined, as is the case of $r_3$ for
small $\lambda_\perp$ in Fig.~\ref{tent_map_rates_trans_lyap1}(c)).

\begin{figure}[htp]
    \centering
   \includegraphics[width=3.4in,height=2in]{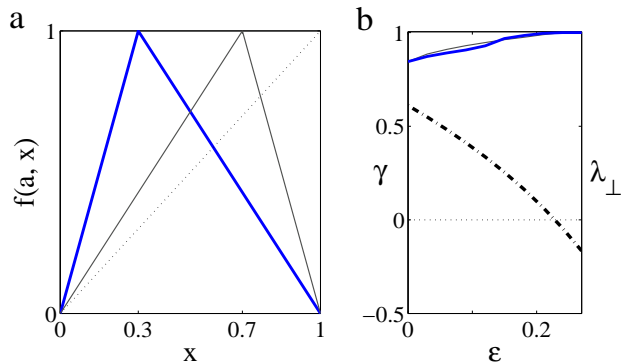}
  \caption{\label{same_lyapunov_tent_map2} (Color online) (a) Two skew tent maps for  $a=0.3$ (thick solid line) and $a=0.7$ (thin solid line), respectively.
  (b) Transversal Lyapunov exponent $\lambda_\perp$ (dash-dotted line) and two synchronization indices $\gamma$ (thick solid line for $a=0.3$ and thin solid line for $a=0.7$). Two $\gamma$ values are almost the same. Transversal Lyapunov exponent $\lambda_\perp$ decreases as $\varepsilon$ increases and crosses 0, when $\gamma$ approaches 1.}
\end{figure}

We now consider some specific parameter values to show that the
transition rates can effectively discriminate two different systems
with identical $\lambda_\perp$ and similar~$\gamma$.
Fig.~\ref{same_lyapunov_tent_map2}(a) shows coupled skew tent maps
for $a=0.3$  and $a=0.7$. Their two Lyapunov exponents are exactly
the same because of their symmetry around $a=1/2$, but the maps are
clearly different and their dynamics are different.
Fig.~\ref{same_lyapunov_tent_map2}(b) illustrates how coinciding
transversal Lyapunov exponents $\lambda_\perp$ (stability property)
and synchronization indices $\gamma$ (average strength of
phase locking) for both maps depend on $\varepsilon$. For these
values of $a$ the threshold value of coupling is $\varepsilon_c
\approx 0.2286$.

\begin{figure}[htp]
    \centering
   \includegraphics[width=3.4in,height=2.8in]{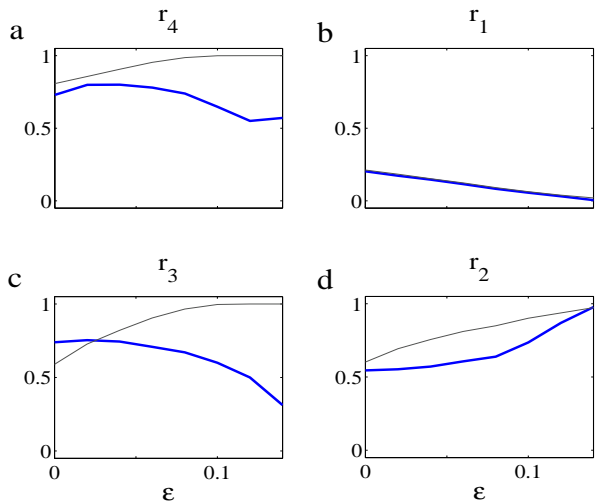}
  \caption{\label{same_lyapunov_tent_map_rate} (Color online) Transition rates for the linearly coupled skew tent maps for $a=0.3$ (thick solid line) and $a=0.7$ (thin solid line). We used 350000 iterations and removed the first 50000 iterations. Differences in the rates $r_{2,3,4}$ are clearly visible.}
\end{figure}

However, Fig.~\ref{same_lyapunov_tent_map_rate} shows that the rates
$r_i$ for $a=0.3$ and $a=0.7$ are different and sometimes
substantially different for a range of $\varepsilon \in (0,
0.15).$ The values of the rate $r_1$ are essentially the same, as
one may expect, because $r_1$ is related to the stability
of synchronized state, which is the same for both cases. The other
rates present a different picture. The system with $a=0.7$ (thin
solid line in Fig.~\ref{same_lyapunov_tent_map_rate}) has higher
$r_{2, 3, 4}$ for almost all range of $\varepsilon $ than those of
$a=0.3$ (thick solid line), which suggests the stronger tendency for
shorter desynchronization events.

As the coupling $\varepsilon$ increases to synchronization threshold
value, the rates $r_{2, 3, 4}$ are undefined while $r_1$ decreases
to zero. This is expected because the system spends most of (if not all)
the time in the vicinity of the synchronized manifold. Thus, the dynamics
is mostly synchronized and the rates are not much relevant. The rate
$r_1 \approx 0$ at $\varepsilon \approx 0.15 < \varepsilon_c \approx
0.2286.$ This probably happens because our partition of the phase space of
the first-return map for the phase difference implies that $r_1=0$
and dynamics is synchronous when the deviation from the complete
synchrony $\theta =0$ is less than $1/2.$ But $\theta(t)$ can still
fluctuate in some smaller limits, which is what happens for
$\varepsilon$ in the interval of about $(0.15, 0.2286).$

This simple case illustrates that although the stability of the
synchronization manifold is the same for two systems, their dynamics
can be substantially different. Lyapunov exponents and
synchronization index $\gamma$ naturally are unable to capture the
differences. However, the analysis using the transition rates can
distinguish two different systems and provide the temporal details
of dynamical behaviors.

\section{\label{intermittency} Analysis of three different intermittency types}

Intermittently synchronized dynamics can be observed near boundary
of synchronization. So far, three types of
intermittencies near the onset of phase synchronization have been
observed: the type-I intermittency~\cite{lee_kwak_lim1}, the eyelet
intermittency~\cite{lee_kwak_lim1, uni_eyelet}, and the ring
intermittency~\cite{ring}.

In this section, we will study these three types of intermittent
behaviors in coupled chaotic oscillators with techniques described
in Section~\ref{methods_sync}. First, we will look at the type-I and
the eyelet intermittencies in bidirectionally coupled R\"{o}ssler oscillators and bidirectionally  coupled
Lorenz oscillators. And then, we will look at the eyelet and the ring intermittencies
in unidirectionally coupled R\"{o}ssler oscillators.

Intermittency types are characterized by the probability
distribution of the laminar (synchronized in our context) phase
length $\ell$ and the scaling behavior of the mean laminar phase
length $<\ell>.$ As far as the type-I and the eyelet intermittencies
are concerned, coupled R\"{o}ssler oscillators and
 coupled Lorenz oscillators with small frequency detuning show the same scaling
properties of synchronous episodes near and away from the
point of bifurcation~\cite{lee_kwak_lim1, uni_eyelet}. However, we will show that there
are important differences in timing of dynamical behavior.

We consider two coupled chaotic R\"{o}ssler oscillators:
\begin{eqnarray}\label{bidirection_Rossler} x^\prime_{1, 2} &=& - \omega_{1, 2} y_{1, 2} - z_{1, 2} +
\varepsilon_{1, 2} (x_{2, 1} - x_{1, 2}), \\ \nonumber %
y^\prime_{1, 2} &=& \omega_{1, 2} x_{1, 2} + 0.15 y_{1, 2}, \\ \nonumber
z^\prime_{1, 2}&=& 0.2 + z_{1, 2} (x_{1, 2} -10), \nonumber \end{eqnarray}
where $\omega_{1, 2}$ are the natural frequencies of each chaotic
oscillator and $\varepsilon_{1, 2}$ measure the strengths of the linear coupling. The phase difference is obtained by
\begin{equation}\label{rossler_phase} \theta = \phi_1 - \phi_2 =
\arctan \left( \frac{y_1}{x_1} \right) - \arctan \left(
\frac{y_2}{x_2} \right).\end{equation}

We also consider two coupled chaotic Lorenz oscillators:
\begin{eqnarray}\label{bidirection_Lorenz}
x^\prime_{1, 2} &=& 10.0(y_{1, 2} - x_{1, 2}) + \varepsilon (x_{2, 1} -
x_{1, 2}), \\ \nonumber %
y^\prime_{1, 2} &=& (36.5 + \gamma_{1, 2}) x_{1, 2} - y_{1, 2} -
x_{1, 2} z_{1, 2}, \\ \nonumber %
z^\prime_{1, 2} &=& -3.0 z_{1, 2} + x_{1, 2} y_{1, 2},\nonumber
\end{eqnarray}
where $\gamma_{1, 2}$ are parameters for the frequency detuning of each chaotic
oscillator and $\varepsilon$ measures the strength of the linear
coupling. As in~\cite{lee_kwak_lim1}, by defining a new variable $u= \sqrt{x^2 + y^2}$ a phase can be defined on the $(u, z)$ plane. The phase difference is obtained by
\begin{eqnarray}\label{Lorenz_phase} \theta &=& \phi_1 - \phi_2 \\
&=& \arctan \left(\frac{z_1-\hat{z}_1}{u_1 - \hat{u}_1} \right)-
\arctan \left( \frac{z_2-\hat{z}_2}{u_2 - \hat{u}_{2}} \right),
\nonumber \end{eqnarray} where $(\hat{u}_{1}, \hat{z}_{1}),
(\hat{u}_{2}, \hat{z}_{2})$ are unstable fixed points in the middle
of the trajectories rotating in the $(u, z)$ plane.

\subsection{\label{type_i_intermittency} Type-I intermittency}

In coupled R\"{o}ssler oscillators,  let us
assume that the difference between $\omega_1$ and $\omega_2$ is
small and the systems are bidirectionally coupled with $\varepsilon_{1,
2} = \varepsilon.$ When $\omega_1 = 1.015$ and $\omega_2 = 0.985,$ Lee
et al.~\cite{lee_kwak_lim1} found the first critical value
$\varepsilon_t = 0.0276$ where $\theta$ increases in a nearly periodic
sequence of $2 \pi$ jumps. They showed that for $\varepsilon <
\varepsilon_t,$ the mean laminar phase length $<\ell>$ obeys the scaling rule
\begin{equation} <\ell>  \sim  (\varepsilon_t - \varepsilon)^{-1/2}. \end{equation}
In coupled Lorenz oscillators, let us again assume that the
difference between $\gamma_1$ and $\gamma_2$ is small. When
$\gamma_1 =1.5$ and $\gamma_2 = -1.5,$ Lee. et
al.~\cite{lee_kwak_lim1} also found the first critical value
$\varepsilon_t=6.7.$ The phase desynchronization occurs with $2 \pi$
phase slip, but it shows $\pm 2 \pi$ irregular phase jumping
behaviors which are different from coupled R\"{o}ssler oscillators.
However, the probability distribution of the laminar phase length
$\ell$ and the scaling rule of the mean laminar phase length
$<\ell>$ are the same as coupled R\"{o}ssler oscillators. The
scaling rule together with the corresponding probability
distributions  of synchronous episodes suggest that these two
systems are in the type-I intermittency.

\begin{figure}[htp]
    \centering
   \includegraphics[width=3in,height=2.8in]{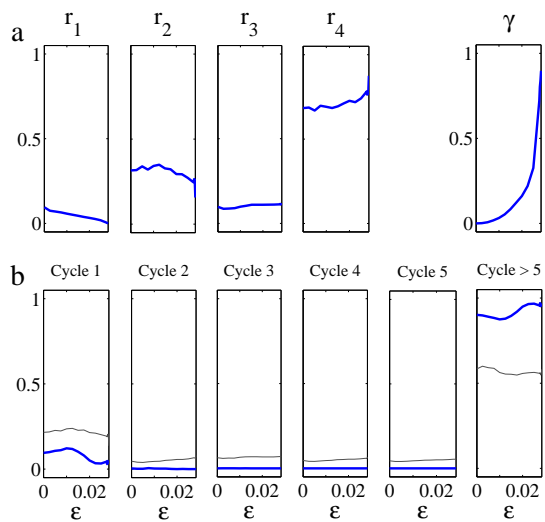}
  \caption{\label{Rossler_Type_I_elelet_best1} (Color online) Synchronized dynamics of bidirectionally coupled R\"{o}ssler oscillators in dependence on the coupling strength 
$\varepsilon$. Two critical values for intermittencies are $\varepsilon_t = 0.0276$ and $\varepsilon_c=0.0286$. (a) Rates $r_i$ and synchronization index $\gamma$. (b) Relative frequencies of the desynchronization events of different durations (durations are measured in cycles of oscillations, i.e. in the number of iterations of the $\phi_{x, i}$ map, see Section~\ref{methods_sync} ). Thick solid line represents the actual relative frequencies of the desynchronization events of different durations, while thin solid line represents the relative frequencies, which would be there under assumption of independent transition rates.  The longer cycles of desynchronization are 
dominant relative to short cycles. }
\end{figure}

\begin{figure}[htp]
    \centering
   \includegraphics[width=3in,height=2.8in]{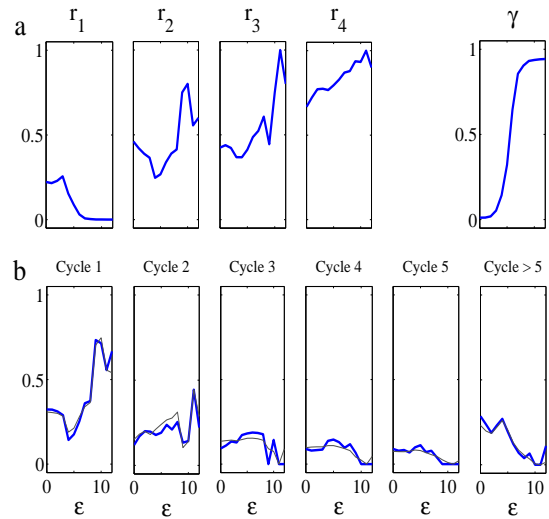}
  \caption{\label{Lorenz_Type_I_elelet_best1} (Color online) Synchronized dynamics of bidirectionally coupled Lorenz oscillators in dependence on the coupling strength $\varepsilon$. Two critical values for intermittencies are $\varepsilon_t =6.7$ and $\varepsilon_c=12$. (a) Rates $r_i$ and synchronization index $\gamma$. The behaviors of the rates and the durations of desynchronization events are different from coupled R\"{o}ssler oscillators.  (b) Relative frequencies of the desynchronization events of different durations (durations are measured in cycles of oscillations, i.e. in the number of iterations of the $\phi_{x, i}$ map, see Section~\ref{methods_sync}). Thick solid line represents the actual relative frequencies of the desynchronization events of different durations, while thin solid line represents the relative frequencies, which would be there under assumption of independent transitions. Unlike Fig.~\ref{Rossler_Type_I_elelet_best1}, short desynchronization events dominate the dynamics. }
\end{figure}

Figs.~\ref{Rossler_Type_I_elelet_best1} and~\ref{Lorenz_Type_I_elelet_best1} show the rates $r_i$ and the
relative frequencies of the desynchronization events of different durations for a range of coupling
strength. By comparing Fig.~\ref{Rossler_Type_I_elelet_best1} with Fig.~\ref{Lorenz_Type_I_elelet_best1}, one can tell that two coupled
systems exhibit different temporal features of the dynamics.

In Fig.~\ref{Rossler_Type_I_elelet_best1}(b), the probability of the long desynchronization events (longer than 5 cycles of oscillations)
is in the range of about $(0.8, 0.9)$ while in Fig.~\ref{Lorenz_Type_I_elelet_best1}(b), this probability is in the
range of about $(0, 0.25).$ This is true, in particular, for similar values of $\gamma$. Since coupled R\"{o}ssler oscillators have high
probability of long cycles, two oscillators spend long time in the desynchronized state once the trajectory leaves synchronous state.
However it does so very rarely. Coupled Lorenz oscillators (Fig.~\ref{Lorenz_Type_I_elelet_best1}) present a different temporal
dynamics. Short desynchronization cycles prevail and the comparable values of $\gamma$ are reached with a larger number of short desynchronized episodes.

The relative frequencies of desynchronization events of different durations can be close to the relative frequencies estimates from the transition rates like Figs.~\ref{Lorenz_Type_I_elelet_best1}(b) and~\ref{uni_eyelet_type_i}(b) or can be different from these estimates like Figs.~\ref{Rossler_Type_I_elelet_best1}(b) and~\ref{ring_intermittency_rates}(b). The transitions between regions are not necessarily independent (or near independent). However, the estimated relative frequencies of desynchronization events of different durations still show some similar tendencies.

\subsection{\label{eyelet_intermittency} Eyelet intermittency}

As the coupling strength increases further from the first critical
value $\varepsilon_t,$ exponentially rare $2 \pi$ phase jumps for
coupled R\"{o}ssler oscillators and irregular $\pm 2 \pi$ phase
jumps for coupled Lorenz oscillators can be observed before the
second critical value $\varepsilon_c$ where the phase
synchronization is established. Since the phase slips are rare, the
scaling rule and probability distributions are different from those
of the type-I intermittency. Lee et al.~\cite{lee_kwak_lim1} found
the second critical values $\varepsilon_c=0.0286$ for coupled
R\"{o}ssler oscillators and $\varepsilon_c = 12$ for coupled Lorenz
oscillators. Note that other parameter values $\omega_{1,2}$ and
$\gamma_{1,2}$ are the same as in the type-I intermittency. For
$\varepsilon_t < \varepsilon < \varepsilon_c,$ the mean laminar
phase length~$<\ell>$ obeys the scaling rule
\begin{equation} \ln(<\ell>) \sim  -(\varepsilon_c - \varepsilon)^{1/2}. \end{equation}
The scaling rule and the corresponding probability distribution of the laminar
phase length suggest both systems are in the eyelet intermittency. Since the eyelet
intermittency features the longer laminar phase than that of  the type-I intermittency, $r_1$ is smaller than that of  the type-I intermittency.

Note that when $\varepsilon > \varepsilon_c,$ the system is in the
first region of the first-return map for phases, that is in the
synchronous region, and the rates $r_{2, 3, 4}$ are undefined while
the rate $r_1 =0.$ However the phase difference may fluctuate in
some smaller limits and synchronization index $\gamma$ may still be
below 1.

Figs.~\ref{Rossler_Type_I_elelet_best1}(a) and~\ref{Lorenz_Type_I_elelet_best1}(a) show that the transition
rates $r_i$ of both coupled  R\"{o}ssler oscillators and coupled Lorenz
oscillators are characteristically distinct from  one another near and
away from the second critical value $\varepsilon_c.$  In fact, the rate $r_3$ for coupled R\"{o}ssler oscillators
is substantially smaller than that of coupled Lorenz oscillators for all range of $\varepsilon.$ This small value of $r_3$ together with 
small value of $r_2 \in (0.25, 0.4)$ for coupled R\"{o}ssler oscillators implies that these two oscillators spend substantially longer time in the desynchronous state
once the trajectory leaves synchronous state than coupled Lorenz oscillators do.
The rates $r_{2, 3, 4}$ for coupled Lorenz oscillators show substantial
variability while those for coupled R\"{o}ssler oscillators do not. Therefore the temporal structure of synchronized/desynchronized events is different in bidirectionally coupled R\"{o}ssler oscillators and Lorenz oscillators even though the overall synchrony can be similar.

Now let us also consider unidirectionally coupled R\"{o}ssler oscillators for the eyelet intermittency. We use Eq.~(\ref{bidirection_Rossler}) with $\varepsilon_1=0$ and $\varepsilon_2 = \varepsilon.$ When control parameters $\omega_1=0.93$ and $\omega_2=0.95,$ Hramov et al.~\cite{uni_eyelet} found  two critical values $\varepsilon_t =0.0345$ and $\varepsilon_c=0.042.$ They showed that when $\varepsilon \in (\varepsilon_t, \varepsilon_c),$ the intermittent behavior of unidirectionally coupled R\"{o}ssler oscillators can be classified as the eyelet intermittency. Note that this intermittent behavior can be also treated as the type-I intermittency with noise~\cite{uni_eyelet}.

\begin{figure}[htp]
    \centering
   \includegraphics[width=3in,height=2.8in]{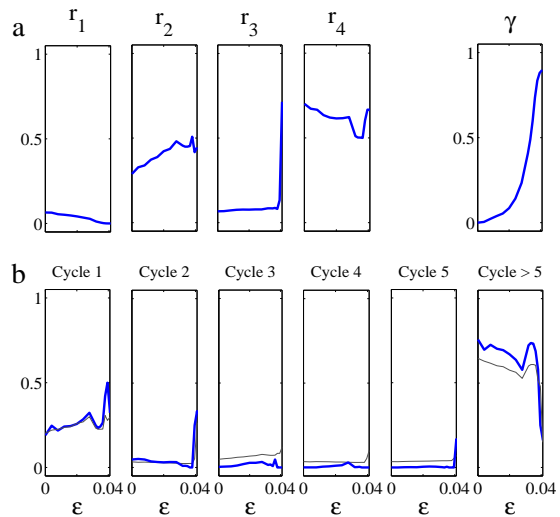}
  \caption{\label{uni_eyelet_type_i} (Color online) Eyelet intermittency in unidirectionally coupled R\"{o}ssler oscillators with parameters $\varepsilon_1=0,$ $\varepsilon_2 =\varepsilon,$ and control parameters $\omega_1=0.93,~\omega_2=0.95$. Two critical points are $\varepsilon_t = 0.0345$ and $\varepsilon_c = 0.042.$ (a) Rates $r_i$ and synchronization index $\gamma$. (b) Relative frequencies of desynchronization events of different durations. Thick solid line stands for actual duration, while thin solid line represents the estimates from the transition rates. The intermittent behavior is different from the bidirectionally coupled R\"{o}ssler oscillators. }
\end{figure}

During the eyelet intermittency, rate $r_{3}$ in Fig.~\ref{uni_eyelet_type_i}(a) experiences substantial variation while that in Fig.~\ref{Rossler_Type_I_elelet_best1}(a) keeps almost the same level. Rate $r_4$ in Fig.~\ref{uni_eyelet_type_i}(a) also shows relatively large variation while that in Fig.~\ref{Rossler_Type_I_elelet_best1}(a) does not. Note that rates $r_1$ for both systems are almost the same. The large variations also can be observed in the durations of desynchronization events for several cycle lengths ($1, 2, >5$) in Fig.~\ref{uni_eyelet_type_i}(b) while almost no variation in Fig.~\ref{Rossler_Type_I_elelet_best1}(b) is observed. Therefore, their temporal structures of synchronized/desynchronized events are different even though the overall synchrony can be similar.

\subsection{\label{ring_intermittency} Ring intermittency}

For the ring intermittency, we consider unidirectionally coupled
R\"{o}ssler oscillators. As in~\cite{ring}, we use
Eq.~(\ref{bidirection_Rossler}) with $\varepsilon_1=0,$
$\varepsilon_2 = \varepsilon,$ and control parameters
$\omega_1=1.0,~\omega_2=0.95.$ The frequency mismatch between two
oscillators is larger than those of the type-I and the eyelet
intermittencies discussed in the previous sections. To explore the
synchronization behavior more clearly, \cite{ring} used the
transformation of coordinates, $\hat{x} = x_2 \cos (\phi_1) + y_2
\sin (\phi_1)$ and $\hat{y} = -x_2 \sin(\phi_1) + y_2 \cos
(\phi_1),$ and plotted the limit cycle in the $(\hat{x}, \hat{y})$
plane.  Then the trajectory on this transformed plane looks like a
ring.

Below the boundary of phase synchronization regime (when
$\varepsilon < \varepsilon_c$), the dynamics of the phase difference
$\theta (t)$ is persistently and intermittently interrupted by
sudden phase slips ($2 \pi$ jumps). As $\varepsilon$ increases
further, the limit cycle starts to enclose the origin near the first
critical point and the phase destruction
 begins. Two critical values $\varepsilon_t = 0.1097 $
and $\varepsilon_c = 0.124$ were observed \cite{ring}, such that (i) the phase synchronization occurs for $\varepsilon
> \varepsilon_c,$ (ii) the (ring) intermittent behavior occurs for
$\varepsilon_t < \varepsilon < \varepsilon_c,$ and (iii) the asynchronous
dynamic occurs for $\varepsilon < \varepsilon_t.$

\begin{figure}[htp]
    \centering
   \includegraphics[width=3in,height=2.8in]{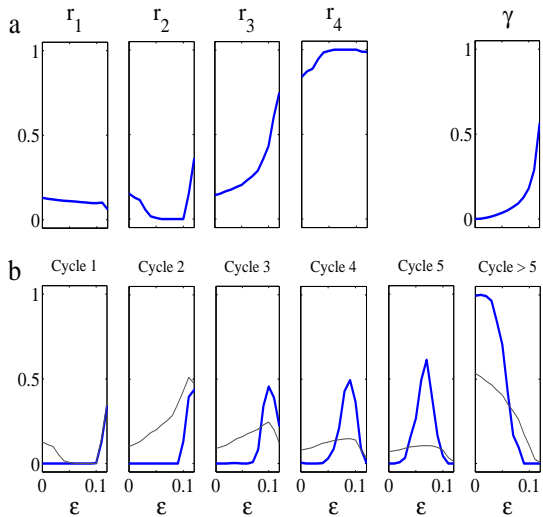}
  \caption{\label{ring_intermittency_rates} (Color online) Ring intermittency in coupled R\"{o}ssler oscillators with parameters $\varepsilon_1=0,$ $\varepsilon_2 =\varepsilon,$ and control parameters $\omega_1=1.0,~\omega_2=0.95$. The ring intermittency occurs between $\varepsilon_t = 0.1097$ and
$\varepsilon_c = 0.124.$ (a) Rates $r_i$ and synchronization index $\gamma$. (b) Relative frequencies of desynchronization events of different durations. Thick solid line stands for actual duration, while thin solid line represents the estimates from the transition rates.}
\end{figure}

Fig.~\ref{ring_intermittency_rates}(a) shows that during the ring
intermittency, the rate $r_1$ slowly decreases. However the rates
$r_{2, 3}$ change drastically during the ring intermittency while
$r_4$ is an almost constant during this parameter regime. These
changes of the rates imply that the probabilities of long cycles are
getting smaller as the coupling strength increases (see
Fig.~\ref{ring_intermittency_rates}(b)). On the other hand, the
probabilities of short desynchronizations (lasting for one or two
cycles of oscillations) are getting higher.

Now let us compare the ring intermittency in Fig.~\ref{ring_intermittency_rates} with the eyelet intermittency in Fig.~\ref{uni_eyelet_type_i}. Both systems are unidirectionally coupled R\"{o}ssler oscillators, but with different control parameters $\omega_{1, 2}$.

Since the phase slips in the eyelet intermittency are rare, as one can expect, rate $r_1$ in Fig.~\ref{uni_eyelet_type_i}(a) has much smaller value than that in Fig.~\ref{ring_intermittency_rates}(a). Rate $r_4$ in Fig.~\ref{uni_eyelet_type_i}(a) has slightly lower value than that in Fig.~\ref{ring_intermittency_rates}(a). Rate $r_2$ in Fig.~\ref{ring_intermittency_rates}(a) shows huge fluctuation while that in Fig.~\ref{uni_eyelet_type_i}(a) does not. The probability of long desynchronization events (Cycle$>5$) for the ring intermittency is close to zero while that in the eyelet intermittency starts to decrease from high level ($>0.5$). Thus, long desynchronization events are dominant for the eyelet intermittency while short desynchronization events are dominant for the ring intermittency.

\section{\label{neuronal_model} Coupled neurons example}

Intermittency and intermittent and otherwise imperfect synchronization can be observed in the brain in different conditions (for example, see~\cite{hurtado2004, park_rubchinsky1, hramov_koronovskii_rats, gong_nikolaev}).
Thus it will be useful to apply the method to more realistic neural-like model systems to study
what rates $r_i$ may reveal about synchronized dynamics. We consider
two mutually coupled neurons described by a single-compartment
conductance-based Hodgkin-Huxley type equations. We couple these
neurons through inhibitory synapses. The details for this kind of a
model can be found in, for example,~\cite{ermentrout_terman}. The
equations for each cell can be written as:
\begin{eqnarray} C_m v^\prime &=& - I_L - I_{Na} - I_K  - I_{syn} + W + I_0,  \\
       x^\prime &=& \varphi (\alpha_x(v) (1-x) - \beta_x (v) x), \nonumber \end{eqnarray}
where $I_L = g_L (v-E_L),$ $I_{Na} = g_{Na} m^3(v) h(v-E_{Na}),$
          $I_K = g_K n^4(v-E_K)$ represent leak, sodium
and potassium currents respectively. Independent
Gaussian white noise $W$ with zero mean and standard deviation of
 $2.12$ is added to both neurons. $I_0$ is an external current. $x$
stands for three different gating variables $m, n$ and $h$ with the
following $\alpha(v)$ and $\beta(v)$:
\begin{eqnarray}
\alpha_n(v)&=& 0.01 (v+55)/(1-\exp(-(v+55)/10)), \\ \nonumber %
\beta_n(v)&=&0.125 \exp(-(v+65)/80), \\  \nonumber %
\alpha_m(v)&=&0.1 (v+40)/(1-\exp(-(v+40)/10)), \\  \nonumber %
\beta_m(v)&=&4 \exp(-(v+65)/18), \\ \nonumber %
\alpha_h(v)&=&0.07 \exp(-(v+65)/20), \\ \nonumber %
\beta_h(v)&=& 1/(1+\exp(-(v+35)/10)). \nonumber  \end{eqnarray}

The term $I_{syn}$ represents the synaptic current. For a cell $i \in \{1, 2\}$, the synaptic current $I_{syn,
i} = g_{syn, i} (v_i - v_{syn}) s_j$ where $j \in \{1, 2\} \setminus \{i\}.$
The synaptic variable $s$ (the fraction of activated synaptic channels) is modeled by the first-order kinetic equation
in the form: \begin{eqnarray}
s^\prime = \alpha (1-s) H_\infty (v -\theta_v) - \beta s,
\end{eqnarray} where $H_\infty(x) = 1/(1+\exp[-x/\sigma_s])$ is a
sigmoidal function. The parameter values are the following: $C_m=1,$
$g_{Na}=120,$ $E_{Na}=50,$ $g_K=36,$ $E_K=-77,$ $g_K=0.3,$
$E_L=-54.4,$ $v_{syn}=-85,$ $I_0=10,$ $\varphi=0.35,$ $\theta_v=0,$
$\sigma_s=5,$ $\alpha=2,$ and $\beta=0.05.$ Both cells are
self-oscillators in the absence of coupling ($g_{syn, i} = 0$).
Voltages $v$ for both neurons exhibit spiky time-series (sequence of
action potentials) and were filtered to the gamma-band ($30-80$~Hz).
The Hilbert transformation was used to define the phases and then
first-return maps were constructed.
\begin{figure}[htp]
    \centering
   \includegraphics[width=3.4in,height=3in]{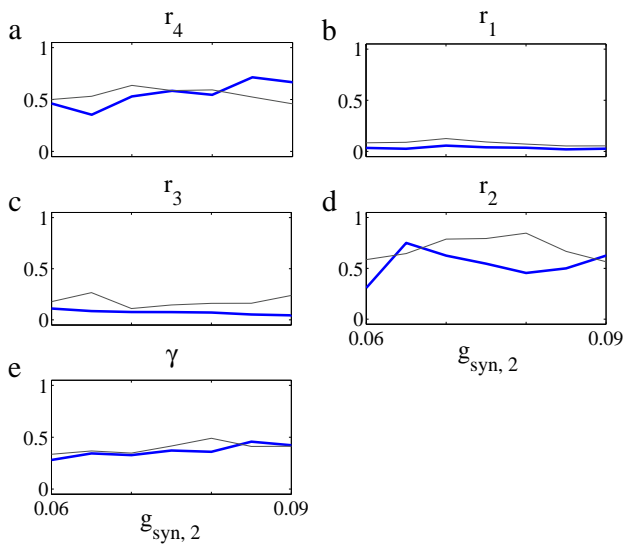}
  \caption{\label{gamma_compare}  (Color online)  Rates $r_i$ and synchronization index $\gamma$ for varied $g_{syn, 2}$ with fixed $g_{syn, 1} = 0.1$. Although the levels of the synchronization index $\gamma$ are similar, the rate $r_{2, 4}$ experience more substantial variation depending on the synaptic strength
  $g_{syn, 2}$. Small noise also exerts some influence on the rates, which is not as prominent in $\gamma$. Thick solid line is the noiseless case and thin solid line is for the noise $W$.}
\end{figure}
 Fig.~\ref{gamma_compare} shows that the synchronization index $\gamma$ changes slowly for a range of values
of $g_{syn, 2}$ (for fixed $g_{syn, 1}=~0.1$), but the rates $r_i$
exhibit more substantial variations depending on $g_{syn, 2}$. As
can be seen from  Fig.~\ref{gamma_compare}(a, d), the rates $r_{2,
4}$ are changing substantially with the similar $\gamma.$ This
implies that the rates describe the fine temporal structure of
synchronization/desynchronization events. The rates (and thus the
fine temporal structure of synchronous/desynchronous dynamics) may
be more affected by noise than synchronization index $\gamma$
(Fig.~\ref{gamma_compare}). Especially for small values of coupling
$g_{syn, 2}$  overall degree of phase locking is not influenced
by small noise. However the rates $r_{2, 3, 4}$ experience
substantial variations and thus corresponding timing of
synchronization/desynchronization events is dependent on the noise
level. Therefore, exploring the fine temporal structures of coupled oscillators in experimental data
may provide some help in constraining multiple alternatives in modeling studies, like in~\cite{park_rubchinsky2}.

\section{\label{discussion} Discussion}

We studied the first-return maps of the phase-difference between two oscillatory signals. If there is some phase locking is present on the average, we can study the phase-locking relationship on each cycle of oscillations. The phase space of this map is partitioned into four regions and the transition rates $r_i$ between regions are computed. These transition rates characterize the details of temporal dynamics, which are missed by standard synchronization measures (such as $\gamma$ in Eq.~(\ref{sync_gamma})).

Using a series of characteristic model systems, such as coupled skew tent maps (for which the values of Lyapunov exponents can be computed analytically), coupled R\"{o}ssler and Lorenz oscillators, and coupled Hodgkin-Huxley neuronal models, we showed that the aforementioned transition rates provide a complimentary description of synchronized dynamics. They describe the dynamics away from the synchronized state, so that even if the stability (transversal Lyapunov exponent) of synchronization manifold and time-series based measure of synchrony (such as a phase-locking index) are the same, the rates will capture the details of desynchronization events. In an extreme case, the same level of synchrony may be reached with a relatively rare but long desynchronization events and numerous short desynchronization events. The rates allow for discrimination of this alternative.

The duration of desynchronization events can be computed from this approach and is related to the rates (although not necessarily defined by them completely). The traditional analysis of intermittent behavior, including intermittency near synchronization onset, considers duration of laminar (i.e. synchronized in the case of synchronization) episodes. These distributions are universal and depend on the type of intermittency. But the same intermittency type in different systems can naturally have different dynamics of synchronization/desynchronization events (as reinjection mechanism is not defined by the type of intermittency). We showed how the rates are able to capture this difference and corresponding difference in the timing of synchronized/desynchronized dynamics.

>From the phase-space-based approach point of view, the system spends a substantial fraction of time away from the synchronization manifold, in other parts of the phase space. Therefore the stability of this manifold and the bifurcations that it can experience provide incomplete information about the dynamics. The introduced rates $r_i$ can help to fill this gap. In the time-series based approach, the global synchronization measures are not designed to describe the relationship between phases on each cycle of oscillations. The approach considered here allows to inspect the phase locking at each cycle of oscillations (naturally, only if some synchronization level is present overall). Recent experimental and modeling results ~\cite{park_rubchinsky1, park_rubchinsky2} indicate that the fine temporal structure of synchronization is important. Here we used standard test systems to explain how the analysis of this structure works.

The methodology considered here can be applied to any appropriate data with some level of phase locking regardless of the mechanisms of synchronous dynamics. We clarified here what this method yields in terms of characterization of the phase space and in terms of the fine temporal structure of phase locking. However, we think in future studies the method maybe used as a tool to study some mechanisms of synchrony or at least distinguish between different mechanisms of synchrony (in particular when parameters of coupling cannot be changed in experiment). 

We would like to conclude the paper with four notes regarding the
method. First, we want to reiterate that for this kind of analysis
the signals should be phase-locked on average. The synchrony is
essentially non-instantaneous phenomenon. Only if there is a preferred phase-locking angle, we can follow deviations from
it on each cycle (and that is yet another reason for necessity of
computations of the global synchronization measures or synchrony
measure over moderate size time-window, like
in~\cite{hurtado_rubchinsky_2004}).

Second, the partitioning of the phase space of the $\phi_{x, i}$ map into four parts is somewhat arbitrary. It simplifies the computation of rates, as there are only four of them. It also implies that if the phases do not deviate from more than $\pi/2$, they are considered in synchrony. This appears to be a reasonable assumption for some systems (for example, smaller tolerance maybe easily destroyed by noise and fluctuations, larger tolerance is less likely to be acceptable for information transmission in the synchronized system). However in certain applications a finer partitioning may be necessary. This is apparently possible to implement, although the rates will be harder to interpret. Moreover, the major idea of the considered approach - to study how the phase difference deviates from the preferred phase-locking angle - can be considered not only in a discrete-time framework as we did here, but also can be generalized to continuous time.

Third, the approach developed here is agnostic to the presence or absence of noise in the time-series. Noise surely can contribute to the generation of imperfect phase locking and affect the transition rates (as in the example in the previous Section). The properties of noise may be reflected in the transition rates and the considered methodology may be a tool to study the effect of noise on synchronous dynamics. This appears to be an important subject for future research.

Fourth, if the system is very close to synchronization threshold, this analysis may not necessarily be very informative, because the system will spend most of the time in the synchronized state. However, as we mentioned in the Introduction, synchronization in living systems may be highly variable. In this case, synchronized events may be relatively rare and/or relatively short (although they may still be significant functionally). Thus the approach considered here (as well as its variations) have  a potential to describe functionally important temporal features of synchronization/desynchronization dynamics.

\begin{acknowledgments}
This study was supported by NIH grant R01NS067200 (NSF/NIH CRCNS).
\end{acknowledgments}

\end{document}